%% file: delayed-cloud.tex
\documentclass[letterpaper, 10 pt, conference]{ieeeconf} 
\IEEEoverridecommandlockouts 
\usepackage{graphicx}
\usepackage{subfigure}
\usepackage{amsmath}
\usepackage{amsfonts}
\usepackage[psamsfonts]{amssymb}
\usepackage{cases}
\usepackage{mathtools}
\usepackage{color}
\usepackage{epstopdf}

\newtheorem{lemma}{Lemma}

\newtheorem{theorem}{Theorem}
\newtheorem{algorithm}{Algorithm}
\newtheorem{assumption}{Assumption}

\newtheorem{problem}{Problem}

\IEEEtriggeratref{5}

\newcommand{\R}{\mathbb{R}}
\newcommand{\N}{\mathbb{N}}

\newcommand{\lve}{L_{\nu, \epsilon}}
\newcommand{\hmve}{\hat{\mu}_{\nu, \epsilon}}
\newcommand{\hxve}{\hat{x}_{\nu, \epsilon}}
\newcommand{\hxvei}{\hat{x}_{\nu, \epsilon, i}}

\newcommand{\ep}{\epsilon}
\newcommand{\dx}{\nabla_{x}}
\newcommand{\dxi}{\nabla_{x_i}}
\newcommand{\dmu}{\nabla_{\mu}}
\newcommand{\pmu}{\Pi_{\mathcal{D}_{\nu}}}
\newcommand{\px}{\Pi_{X}}

\newcommand{\dv}{\mathcal{D}_{\nu}}
\newcommand{\bx}{\bar{x}}

\newcommand{\hmvej}{\hat{\mu}_{\nu, \ep, j}}
\newcommand{\grad}{\nabla}

\newcommand{\ns}{\!}

\newcommand{\qve}{q_{\nu, \epsilon}}

\title{\LARGE \bf
Cloud-Based Centralized/Decentralized Multi-Agent Optimization \\ with Communication Delays
}

\author{Matthew T. Hale,$^\dag$ Angelia Nedi\'{c},$^\star$ and Magnus Egerstedt$^\dag$
\thanks{$^\dag$School of Electrical and Computer Engineering, Georgia Institute of
Technology, Atlanta, GA 30332, USA. Email: \texttt{\{matthale, magnus\}@gatech.edu}. Research supported in
part by  the NSF under Grant CNS-1239225.}
\thanks{$^\star$Department of Industrial and Enterprise Systems Engineering, University of Illinois, Urbana IL, 61801, USA.
Email: \texttt{angelia@illinois.edu}.}
}

\begin{document}
\maketitle
\thispagestyle{empty}
\pagestyle{empty}

\input{0-abstract.tex}
\input{1-introduction.tex}
\input{2-problem-an.tex}

\input{3-architecture-an.tex}

\input{4-convergence-an.tex}

\input{5-results.tex}

\input{6-conclusion.tex}

\bibliographystyle{plain}{}
\bibliography{sources}

\end{document}

%% file: 0-abstract.tex
\begin{abstract}
We present and analyze a computational hybrid architecture for performing multi-agent optimization.
The optimization problems under consideration have convex objective and constraint
functions with mild smoothness conditions imposed on them. 
For such problems, we provide a primal-dual algorithm implemented in the hybrid architecture,
which consists of a decentralized network of
agents into which centralized information is occasionally injected, and we establish its convergence properties. 
To accomplish this, a central cloud computer
aggregates global information, carries out computations of the dual variables based on this information, 
and then distributes the updated dual variables to the agents. 
The agents update their (primal) state variables and also communicate  
among themselves with each agent sharing and receiving state information with
some number of its neighbors. Throughout, 
communications with the cloud are not assumed to be synchronous or instantaneous, and 
communication delays are explicitly accounted for in the modeling and analysis of the system. 
Experimental results are presented to support the theoretical developments made. 
\end{abstract}

%% file: 1-introduction.tex
\section{Introduction}
Algorithms for multi-agent and distributed 
optimization have been considered for a variety of problem formulations
in part because of the varied collection of application domains in which
such problems arise. Applications of multi-agent optimization can be found in robotics~\cite{guo02,cortes05,soltero13},
power systems~\cite{nazari14}, sensor networks~\cite{trigoni12, khan09, rabbat04}, 
and communications~\cite{chiang07, wei13, kelly98, mitra94}. 

These diverse applications lead to optimization problems of many 
different formulations. Correspondingly, algorithms
have been developed that allow for a broad range of problem characteristics.
For example, in \cite{notarstefano09} 
distributed linear programs with constraints on network connectivity and memory
are considered. 
In \cite{gharesifard14} the authors consider a distributed method for
minimizing a sum of convex functions over a digraph.
The authors of \cite{nedic09} devise an algorithm for minimizing
a sum of possibly non-differentiable convex functions over a time-varying
digraph. Problems with time-varying graphs, non-differentiable objective functions, and
noisy communication links are considered in 
\cite{srivastava11}, while \cite{lobel11} considers similar problems in which nodes
in the network are expected to fail over time. 

The development of decentralized methods in
optimization has been motivated in part by the fact that centralized methods
may not scale well for very large networks of agents~\cite[Section 1.1]{tsitsiklis84}. At the same time,
centralized methods can more efficiently solve some problems, like the Credit Assignment Problem
in multi-agent robotics, than decentralized methods~\cite[Section 3.1]{panait05}. 
Such examples indicate that
centralized information may be rich in a way that could be useful in networks
which would otherwise be purely decentralized. 
In adding a centralized component to a decentralized network, it seems likely that the centralized component would
operate slower than the decentralized components. Nevertheless, 
one may ask whether it would be useful to occasionally inject global information into a multi-agent network where
such information would otherwise be absent. 

Towards answering this question, 
we present here a multi-agent optimization architecture in which a cloud computer is used
to occasionally provide centralized information 
to a network of agents solving a nonlinear programming problem.
This cloud-based optimization architecture was
introduced in~\cite{hale14}, though
here we substantially 
broaden the class of problems to be solved and allow for communications delays 
when communicating with the cloud.
The cloud carries out computations based on information sent to it
by agents in the network, and intermittently disseminates these results to the agents for use in their
own computations. At the same time, each agent shares its state with some number of neighboring
agents at each time. 
In~\cite{hale14}, the assumption was made that all information in the network was synchronized
at each time, so that all computations were relying on the same information.
Here we eliminate this assumption and, as a consequence, delays occur which give rise to
various kinds of errors. We present explicit bounds on these errors in terms of network constants and algorithm parameters.

To solve nonlinear program in a distributed manner  across many agents, we
cast such problems as a variational inequality and then use Tikhonov regularization, 
 which endows the resulting variational inequality with certain properties that let us draw
from existing convergence results. 
In particular we consider a fixed regularization as was done in \cite{koshal11}.
A fixed regularization is desirable for multi-agent problems because it may be difficult to
synchronize the timing of the changes in regularization parameters across large networks. 
Accordingly, we use the approach of \cite{koshal11} as a starting point, though the problem
and approach there are quite different from the current paper. 
In~\cite{koshal11}, the need for these results stems from reducing computation times
in Lagrangian subproblems associated with a dual optimization scheme
by allowing inexactness in some computations.
Here we use a different architecture and different model for delays
to operate as fast as possible by using
the most recent information available to the agents. Doing so
results in delays in communication, and  computations that rely
on information of different ages. The resulting structure of delays
will be detailed below. 
  
The rest of the paper is organized as follows. 
Section~\ref{sec:global}  will cover the background concerning the
problem to be solved and a centralized method for solving it.
Then, Section~\ref{sec:arch} will cover the cloud architecture and modify the centralized solution method
to fit with the hybrid centralized/decentralized system.
Section~\ref{sec:analysis} will present the convergence results and error bounds of the partially decentralized algorithm, 
and Section~\ref{sec:numer} will present experimental results from an implementation of this algorithm on a team of mobile robots. 
Section~\ref{sec:concl} will conclude the paper. 

%% file: 2-problem-an.tex
\section{Problem Formulation and \\ Centralized Solution}\label{sec:global}
In this section we formulate the problem to be solved.
This section treats everything globally and the results here will be
modified later in Section \ref{sec:arch} to fit with a hybrid
architecture described therein. 

\subsection{Variational Inequality Setup} \label{ss:probsetup}
Consider a multi-agent optimization problem comprised of $N$ agents indexed
by ${i \in I := \{1, \ldots, N\}}$. Suppose that agent $i$ has state $x_i \in \R^{n_i}$
with $n_i \in \N$, and
let the vector $x$ denote the column vector of all states, namely
\begin{equation*} 
x = \left(\begin{array}{c} x_1 \\ x_2 \\ \vdots \\ x_n \end{array}\right) \in \R^n,
\end{equation*}
where $n = \sum_{i=1}^{N} n_i$. 
Let each agent have a local objective function defined only on its own state, $f_i(x_i)$. We assume that 
$f_i : \R^{n_i} \to \R$ is $C^1$ and convex. We also consider a global cost
that is not necessarily separable, $c(x)$,
and assume that $c : \R^n \to \R$ is both $C^1$ and convex as well. 
We assume that agent $i$ knows
$c$ and $f_i$, but not $f_{\ell}$ for any $\ell \neq i$; that is, each agent
knows the non-separable cost and its own local cost, but not the local cost function of any other agent. 
We assume further that the cloud does not know $f_i$ for any $i \in I$. 

The agents are collectively subject
to global inequality constraints of the form
\begin{equation*}
g(x) = \left(\begin{array}{c} g_1(x) \\ g_2(x) \\ \vdots \\ g_m(x) \end{array}\right) \leq 0,
\end{equation*}
where $g : \R^n \to \R^m$ with $m\ge 1$. The constraint functions
$g_j : \R^n \to \R$ are assumed to be convex and $C^1$ for all $j \in J := \{1, \ldots, m\}$. 
In addition to $g$,
each agent's state is also constrained to lie in a given set $X_i\subset\R^{n_i}$, i.e.,
\begin{equation*}
x_i \in X_i
\end{equation*}
for every $i$, where $X_i$ is non-empty, compact, and convex. Letting
\begin{equation*}
X = X_1 \times X_2 \times \cdots \times X_N,
\end{equation*}
we encapsulate each set constraint by requiring
\begin{equation*}
x \in X. 
\end{equation*}
Regarding the class of optimization problems under consideration, we
summarize the conditions that we have imposed in the following assumption.
\begin{assumption} \label{as:functions}
The set $X$ is non-empty, compact, and convex. The functions $\{g_j\}_{j \in J}$
and $c$ are convex and $C^1$ in $x$, and $f_i$ is convex and $C^1$ in $x_i$ for all $i \in I$. 
\hfill $\blacktriangle$
\end{assumption}

For notational convenience,
define the function 
\begin{equation*} 
f(x) = \sum_{i=1}^{N} f_i(x_i) + c(x).
\end{equation*}
Let $\nabla_{x_i}$ denote the operator $\frac{\partial }{\partial x_i}$ and define the map
\begin{equation*}
\nabla f(x) = \Big(\nabla_{x_1}\big(f_1(x_1) + c(x)\big), \ldots, \nabla_{x_N}\big(f_N(x_N) + c(x)\big)\Big).
\end{equation*}
We enforce the following assumption on $\nabla f$.
\begin{assumption} \label{as:gradlip}
The map $\nabla f$ is Lipschitz continuous with Lipschitz constant $L_f$. \hfill $\blacktriangle$
\end{assumption}

Note that any collection of functions $f_i$ and $c$ which are all $C^2$ 
comprise an $f$ that automatically satisfies
Assumption~\ref{as:gradlip} whenever $X$ is compact 
(cf.\ Assumption~\ref{as:functions}).
Concerning the constraints $g$, we have the following assumptions. 
\begin{assumption}\label{as:slater} (Slater's Condition)
There exists a vector $\bar{x} \in X$ such that $g(\bar{x}) < 0$, i.e., the constraints
are strictly feasible at $\bar{x}$. \hfill $\blacktriangle$
\end{assumption}
\begin{assumption} \label{as:gradglip}
The gradient of each constraint, $\nabla g_j$, $j \in J$, is Lipschitz with constant $L_j$ and hence
$\nabla g$ is Lipschitz with constant 
\begin{equation*}
L_g = \sqrt{\sum\limits_{j=1}^{m} L_j^2}.
\end{equation*}
\hfill $\blacktriangle$
\end{assumption}

A global formulation of the multi-agent optimization problem under consideration is given by:
\begin{problem} \label{prob:global}
\begin{equation*} 
\begin{aligned} 
\textnormal{minimize } f(x)& \\
\textnormal{subject to } g(x) & \leq 0 \\ 
                            x  \in &X.
\end{aligned}
\end{equation*}
\hfill $\Diamond$
\end{problem} 

Assumption~\ref{as:functions} guarantees that Problem~\ref{prob:global} has a solution and
Assumption~\ref{as:slater} guarantees that a dual solution exists with no duality gap. 
Denote an optimal primal-dual pair for Problem~\ref{prob:global} by $(\hat{x}, \hat{\mu})$. 
We now define the Lagrangian associated with Problem~\ref{prob:global} as
\begin{equation*}
L(x, \mu) = f(x) + \mu^Tg(x),
\end{equation*}
where $x \in X$ is as defined above and $\mu$ is a vector of
Kuhn-Tucker multipliers in the non-negative orthant of $\R^m$,
denoted $\R^m_{+}$.
By definition, $L(\cdot, \mu)$ is convex for all $\mu \in \R^m_{+}$ and $L(x, \cdot)$ is concave for every $x$. 
Seminal work by Kuhn and Tucker \cite{kuhn51} showed
that optimal primal-dual pairs for Problem \ref{prob:global} are saddle points of $L$. This
saddle point condition can be expressed concisely as: for all $x \in X$  
and $\mu \in \R^m_{+}$,
\begin{equation*} 
L(\hat{x}, \mu) \leq L(\hat{x}, \hat{\mu}) \leq L(x, \hat{\mu}).
\end{equation*}

The problem of finding Lagrangian saddle points can be restated as a variational inequality (e.g., \cite[Section 11.1]{konnov07}). 
Let $\nabla_{x}$ and $\nabla_{\mu}$ denote the operators $\frac{\partial}{\partial x}$ and
$\frac{\partial}{\partial \mu}$, respectively. 
In the variational inequality setting, we wish to find a point $(\hat{x}, \hat{\mu}) \in X \times \R^m_{+}$ such that
\begin{equation*}
\left[\left(\begin{array}{c} x \\ \mu \end{array}\right) - \left(\begin{array}{c} \hat{x} \\ \hat{\mu} \end{array}\right)\right]^T
\left[\begin{array}{r} \nabla_{x} L(\hat{x}, \hat{\mu}) \\[3pt] -\nabla_{\mu}L(\hat{x}, \hat{\mu}) \end{array}\right] \geq 0
\end{equation*} 
for all $(x, \mu) \in X \times \R^m_{+}$. 
In order to make use of certain established results concerning variational inequalities, we take two 
further theoretical steps: first we modify the map $(\nabla_x L \,\, -\nabla_{\mu} L)^T$ to make it strongly monotone, 
and second we find a (non-empty) compact, convex set containing the optimal primal-dual vectors. 
We describe both steps below. 

\subsection{Tikhonov Regularization}
The gradient map
\begin{equation*}
\left(\begin{array}{r} \nabla_{x} L(x, \mu) \\[3pt] -\nabla_{\mu}L(x, \mu) \end{array}\right)
\end{equation*}
is monotone, and we use a fixed Tikhonov regularization in order to work with a strongly monotone map. 
This is done by regularizing the Lagrangian function as follows:
\begin{equation*}
\lve(x, \mu) = f(x) + \frac{\nu}{2}\|x\|^2 + \mu^Tg(x) - \frac{\epsilon}{2}\|\mu\|^2,
\end{equation*}
where $\nu > 0$ and $\epsilon > 0$, and these values are kept fixed to avoid the need to synchronize changes in
parameter values across many agents. 

Under this regularization, we see that $\lve(\cdot, \mu)$ is strongly convex for all $\mu \in \R^m_{+}$ and $\lve(x, \cdot)$ is strongly concave for all $x$. 
These properties imply that $\nabla_{x} \lve$ and $-\nabla_{\mu} \lve$ are both strongly monotone; hence,
the map $(\nabla_{x} \lve \,\, -\nabla_{\mu} \lve)^T$ is also strongly monotone. 
In addition, the strongly convex-strongly concave property of $\lve$, together
with Assumption~\ref{as:functions} and Assumption~\ref{as:slater},
guarantee the existence of 
 a unique optimal primal-dual pair,
$(\hxve, \hmve) \in X \times \R^m_{+}$. 
Using the regularized Lagrangian, 
we now state the variational inequality
of interest.
\begin{problem}\label{prob:vi} 
Find the point $(\hxve, \hmve) \in X \times \R^m_{+}$ such that
for all $(x, \mu) \in X \times \R^m_{+}$,
\begin{align*}
\left[\left(\begin{array}{c} x \\ \mu \end{array}\right) - \left(\begin{array}{c} \hxve \\ \hmve \end{array}\right)\right]^T
\left[\begin{array}{r} \nabla_{x} \lve(\hxve, \hmve) \\[3pt] -\nabla_{\mu}\lve(\hxve, \hmve) \end{array}\right] \geq 0.
\end{align*} 
\hfill $\Diamond$
\end{problem}
We note that determining the solution to the above variational inequality is the same as 
determining the saddle-point of the regularized Lagrangian function $\lve$, i.e.,
$(\hxve, \hmve) \in X \times \R^m_{+}$ solves the above variational inequality problem if and only if 
for all $(x,\mu)\in X \times \R^m_{+}$,
\begin{align}\label{eq:saddle}
\lve(\hxve, \mu)\le \lve(\hxve, \hmve) \le \lve(x, \hmve).
\end{align}

\subsection{Bounds on Dual Variables}
We proceed along the lines of \cite{uzawa58} and derive a bound on $\hmve$.
Letting $\bar{x}$ denote a Slater
point for the constraints $g$, 
we define the dual function associated with $\lve$ as 
\begin{equation*}
\qve(\mu) := \min_{x \in X} \lve(x, \mu).
\end{equation*}
Now consider an arbitrary multiplier $\tilde{\mu} \in \R^m_{+}$
and let the point $\tilde{x} \in X$ be such that
$\tilde{x} = \min_{x \in X} \lve (x, \tilde{\mu})$. By the definition of $\qve$ we then have
\begin{equation*}
\qve(\tilde{\mu}) = \lve(\tilde{x}, \tilde{\mu}) \leq \lve(\hxve, \tilde{\mu}). 
\end{equation*}
In view of the saddle-point property of $(\hxve, \hmve)$ (cf.\ Equation~\eqref{eq:saddle}),
it further follows that
\begin{equation*}
\qve(\tilde{\mu}) 
\leq \lve(\hxve, \hmve) \leq \lve(\bar{x}, \hmve).
\end{equation*}
Using the regularized Lagrangian expression, we have
\begin{align*}
\qve(\tilde{\mu}) &\leq f(\bar{x}) + \hmve^Tg(\bar{x}) + \frac{\nu}{2}\|\bar{x}\|^2 - \frac{\epsilon}{2}\|\hmve\|^2 \\
     &\leq f(\bar{x}) + \hmve^Tg(\bar{x}) + \frac{\nu}{2}\|\bar{x}\|^2.
\end{align*}
Rearranging terms then gives
\begin{equation*}
-\sum_{j=1}^{m} \hmvej g_j(\bar{x}) \leq f(\bar{x}) + \frac{\nu}{2}\|\bar{x}\|^2 - \qve(\tilde{\mu}), 
\end{equation*}
from which we conclude that for any $\tilde{\mu} \in \R^m_{+}$,
\begin{equation} \label{eq:firstmubound}
\sum_{j=1}^{m} \hmvej \leq \frac{f(\bar{x}) 
+ \frac{\nu}{2}\|\bar{x}\|^2 - \qve(\tilde{\mu})}{\min\limits_{1 \leq j \leq m} \left\{-g_j(\bar{x})\right\}}.
\end{equation}
For any $\nu > 0$ and $\epsilon > 0$, we certainly have
\begin{equation*}
\qve(\tilde{\mu}) = \lve(\tilde{x}, \tilde{\mu}) \geq L_{0, \epsilon}(\tilde{x}, \tilde{\mu}) = f(\tilde{x}) + \tilde{\mu}^Tg(\tilde{x}) - \frac{\epsilon}{2}\|\tilde{\mu}\|^2.
\end{equation*}
Selecting $\tilde{\mu} = 0$, observe that
\begin{equation*}
L_{0, \epsilon}(\tilde{x}, 0) = f(\tilde{x}).
\end{equation*}
Letting $f_{X}^* = \min_{x \in X} f(x)$, the bound in Equation \eqref{eq:firstmubound} can be simplified to
\begin{equation*}
\sum_{j=1}^{m} \hmvej \leq \frac{f(\bar{x}) + \frac{\nu}{2}\|\bar{x}\|^2 - f_X^*}{\min\limits_{1 \leq j \leq m} \left\{-g_j(\bar{x})\right\}}.
\end{equation*}

Using that $\hmvej \geq 0$ for all $j \in J$ and defining
\begin{equation*}
\dv = \left\{\mu \in \R^m_+ : \|\mu\|_{1} \leq \frac{f(\bx) + \frac{\nu}{2}\|\bx\|^2 
- f_X^* }{\min\limits_{1 \leq j \leq m} \left\{-g_j(\bx)\right\}}\right\},
\end{equation*}
we see that $\dv$ is non-empty, compact, and convex, and we are guaranteed that $\hmve \in \dv$. 
Using $\dv$, we can now define the algorithm used to solve Problem~\ref{prob:vi}.
\begin{algorithm} \label{alg:global}
Given an initial point $(x(0), \mu(0)) \in X \times \dv$, execute the update law
\begin{equation} \label{eq:globalgx}
x(k+1) = \px\Big[x(k) - \alpha \dx\lve(x(k), \mu(k))\Big] 
\end{equation}
\begin{equation} \label{eq:globalgmu}
\mu(k+1) = \pmu\Big[\mu(k) + \tau \dmu\lve(x(k), \mu(k))\Big],
\end{equation}
until some stopping criterion is reached. \hfill $\blacklozenge$
\end{algorithm}

Here $\px[\cdot]$ and $\pmu[\cdot]$ are the projections onto the sets $X$ and $\dv$, respectively, with respect to the standard Euclidean norm.
In the next section we will explicitly
reformulate Algorithm~\ref{alg:global} for the cloud-based multi-agent case.

%% file: 3-architecture-an.tex
\section{Cloud Architecture and Hybrid Solution} \label{sec:arch}
We now cover the architecture that will be used to implement a hybrid
form of Algorithm \ref{alg:global}. Then we cover how the centralized 
means of solving this problem will be modified
and deployed on this architecture. 

\subsection{Overview}
Let agent $i$ have neighborhood set $N_i$ containing the indices of all agents it is directly coupled to
by $g$ and $c$. 
That is, if the computation of
\begin{equation*}
\frac{\partial \lve}{\partial x_i}
\end{equation*}
requires $x_j$, then $j \in N_i$ and agent $j$ sends its state to agent $i$ at each time. 
This structure of communications necessitates that $i \in N_j$ if and only if $j \in N_i$. 

Let the agents share their states with each other at each timestep via 
wireless communication links.
In this framework, each agent stores and manipulates a local copy of 
Problem \ref{prob:vi} onboard and updates its own state within that local copy based on computations
it performs onboard. Within each timestep, agent $i$ computes an updated value of its own state, $x_i$,
and then shares the new value of $x_i$ with agent $j$ for all $j \in N_i$. 
In many cases, including in very large networks of agents, we expect that the neighborhood
set of each agent will be a small subset of total collection of agents so that
$|N_i| < N$ for all $i$, where $|\cdot|$ denotes cardinality. 
Computing values of dual variables using Algorithm
\ref{alg:global} 
will require all states in the system (see Equation \eqref{eq:globalgmu}),
and given that $|N_i| < N$, we see  
that no agent will be able to perform these computations.
Furthermore,
there is no assumption that the communication graph of the system is connected,
nor is it even assumed that each agent is coupled to any other agent. In such cases, there is not any way
to aggregate all states in the network onboard a single agent, even after long periods of time. 

To fill this gap,
we use a cloud computer
as was done in \cite{hale14}. The cloud computer is assumed to be 
capable of executing computationally intensive calculations quickly as would be the case
with a computer cluster or server farm.
Occasionally every
agent sends its state to the cloud and after some time each agent receives back an updated dual vector
that it stores onboard and incorporates into its own local calculations of state updates. 
Because the cloud must aggregate all states in the network, 
it is assumed that there are
delays in communicating with the cloud
and, due to these delays, Algorithm \ref{alg:global} will be modified.

The precise update law used by each agent will be detailed below, though for the current discussion
it is sufficient to note that each agent executes some onboard update law using its own state information,
information it receives from its neighbors, and the most recent dual vector it has received from the cloud, regardless of how long
ago that dual vector was received. 
Suppose the agents send their states to the cloud at some timestep $k_0$ and
suppose they all have some dual vector $\mu_0$ onboard which was received just prior to sending their states to the cloud 
(the states sent at time $k_0$ were not computed using $\mu_0$). 
While the agents are waiting to receive an updated dual vector, $\mu_1$, they continue to communicate with
each other as before and continue to use $\mu_0$, which is held constant onboard
each agent, in their computations of state updates. 

Suppose the agents' states from $k_0$ arrive at the 
cloud\footnote{The agents' states can arrive at the cloud at different times in which
case the cloud can wait to compute
an updated dual vector until it has received all agents' states. Here, we assume all states arrive simultaneously
for simplicity.} 
at time $k_0 + p_0$ for some $p_0 \in \N$. With all states received, the cloud can compute the next
dual update using a rule similar to that in Algorithm \ref{alg:global}. Suppose that computing the next dual vector
takes some number of timesteps $q_0 \in \N$ so that $\mu_1$ has been computed at time $k_0 + p_0 + q_0$. 
Then the cloud sends $\mu_1$ to each agent and it takes $r_0 \in \N$ timesteps to reach the agents,
arriving at time $k_0 + p_0 + q_0 + r_0$. 
Before agent $i$ computes a primal update of $x_i$ using $\mu_1$, it again
sends its state to the cloud and then uses $\mu_1$ in its subsequent computations. Then this process of the agents 
sharing states with their neighbors, 
receiving a delayed dual update, and sending their states to the cloud
 is repeated. Note that the delays $p_0$, $q_0$, and $r_0$
are not assumed to be constant but instead are associated with $k_0$ and are allowed
to vary with each communications cycle, i.e., if the agents again send their states to the 
cloud at time $k_1$, there is no need for $p_1 = p_0$, $q_1 = q_0$, or $r_1 = r_0$. 

\subsection{Update Law Derivation}
To derive the per-agent update law based on the above communications scheme,
we consider the transfer of information through the network starting from the initial point $(x(0), \mu(0))$
and generalize it to an arbitrary communications cycle.

Define the delay between
the agents receiving $\mu(t)$ and $\mu(t+1)$
as $d(t) := p_t + q_t + r_t$.
When the system is initialized, agent $i$ has onboard its own state, $x_i(0)$, 
the primal stepsize $\alpha$, and the initial dual vector, $\mu(0)$. The cloud is initialized with the same value of $\mu(0)$ and stepsize $\tau$. 
At time $k = 0$, each agent sends its state to the cloud and then the optimization
process begins. Upon receiving all agents' states (which comprise the full vector $x(0)$), the cloud will compute $\mu(1)$
according to
\begin{equation*}
\mu(1) = \Pi_{\dv}\big[\mu(0) + \tau \nabla_{\mu} \lve(x(0), \mu(0))\big].
\end{equation*}
Simultaneously, 
agent $i$ will use the update law\footnote{For consistency of notation, $\nabla_{x_i} \lve$ is written as an argument
of the full state vector $x^0(k)$ here, though $\nabla_{x_i} \lve$ will only depend upon the states of agents
with indices in $N_i$.}
\begin{equation} \label{eq:schedx}
x_i(k+1) = \Pi_{X_i}\big[x_i(k) - \nabla_{x_i} \lve(x(k), \mu(0))\big]
\end{equation}
until it receives $\mu(1)$ from the cloud, and updates of this form occur synchronously across the agents. 
When each agent receives $\mu(1)$, it sends the state
$x_i(d(0))$ to the cloud before performing any computations
with $\mu(1)$. After every state of the form $x_i(d(0))$ arrives at the cloud (which is not assumed to happen
instantly), the cloud will compute $\mu(2)$ according to 
\begin{equation} \label{eq:schedmu}
\mu(2) = \Pi_{\dv}\bigg[\mu(1) + \tau \nabla_{\mu} \lve(x(d(0)), \mu(1))\bigg].
\end{equation}

Equation \eqref{eq:schedx} shows that Equation \eqref{eq:globalgx} in Algorithm
\ref{alg:global} can be distributed among the agents in a natural way by having each agent
update its own state (recall that $X=X_1\times\cdots\times X_N$). 
 However, Equation \eqref{eq:schedmu} reveals that
Equation \eqref{eq:globalgmu} must be modified to account for the communication delays present
in the system. In particular, the arguments of $\nabla_{\mu} \lve$ are no longer aligned
in time because $\mu(2)$ is computed using $\mu(1)$ and values of $x$ that are based on $\mu(0)$. 

\subsection{Multi-agent Implementation}
To compactly express this cycle of communication and computation, we implement a change in notation.
Let $\mu$ be indexed by the time variable $t \in \N$.
In this notation, 
the state of each agent will be indexed both over timesteps at which the agents compute primal updates
 and also over which dual vector is currently
being used in its computations. The time index of the agents' computations will be $k \in \N$ and agent 
$i$'s state will have a superscript to denote the time index of the dual variable agent $i$ currently
has onboard. The results of agent $i$'s $k^{th}$ state update using $\mu(t)$ 
will be denoted $x_i^t(k)$. Using this notation, we 
restate Algorithm \ref{alg:global} to explicitly specify the update law for agent $i$ and
to account for the delays in dual vectors seen above. 
\begin{algorithm} \label{alg:dist}
Let agent $i$ have initial state $x^0_i(0)$, stepsize $\alpha$, and initial
dual vector $\mu(0)$. Let the cloud have initial 
multiplier vector $\mu(0)$ and stepsize $\tau$. Execute for each agent $i$, the following two steps:
for all $k = 0, \ldots, d(t) -1, \,\, t \geq 0$, 
\begin{align}\label{eq:distalgx}
x_i^{t}(k+1) &= \Pi_{X_i}\Big[x_i^t(k) - \alpha\dxi\lve\big(x^t(k), \mu(t)\big)\Big] 
\end{align}
and for all $t \geq 0$,
\begin{align}\label{eq:distalgmu}
\mu(t+1) =
\pmu\Big[\mu(t) + \tau\dmu\lve\big(x^{t-1}\big(d(t-1)\big), \mu(t)\big)\Big] 
\end{align}
until a certain stopping criterion is reached by each agent and the cloud. \hfill $\blacklozenge$
\end{algorithm}
In the setting of Algorithm~\ref{alg:dist}, we define $x^{-1}\big(d(-1)\big) = x^0(0)$.

%% file: 4-convergence-an.tex
\section{Convergence Analysis} \label{sec:analysis}
We now show that Algorithm \ref{alg:dist} ``nearly'' converges 
to the saddle point of $\lve$ and bound the quantities
$\|x^t(k) - \hxve\|$ and $\|\mu(t) - \hmve\|$, 
where $x^t(k)$ is obtained by stacking the vectors $x_i^t(k)$, $i=1,\ldots,N$. 
To do so we first establish the following lemma concerning
convergence when $\mu(t)$ is fixed. 

\begin{lemma} \label{lem:geoconv}
Let Assumptions 1-4 hold. 
Let each agent use the primal regularization parameter $\nu$.
Define the constants $M_\nu = \max_{\mu \in \dv} \|\mu\|$ and
$C_f = L_f + \nu + M_{\nu} L_g$. 
Enforce that the primal stepsize $\alpha$ satisfies $0 < \alpha < 2/C_f$. 
Let 
$x^t_{*, i}$ denote the fixed point of Equation \eqref{eq:distalgx} for agent $i$ with $\mu(t)$ fixed,
and define $x^t_* = (x^{t,T}_{*,1}, \ldots, x^{t,T}_{*,N})^T$. 
Define also $x^t(k) = (x^t_1(k)^T, \ldots, x^t_N(k)^T)^T$. 
For a fixed dual vector $\mu(t)$ the sequence $\{x^t(k)\}$, $k = 1, \ldots, d(t)$ satisfies
\begin{equation*}
\|x^t(k) - x^t_{*}\| \leq q_p^{k/2} \|x^t(0) - x^t_{*}\|,
\end{equation*}
where $(0, 1) \ni q_p := 1 - \alpha\nu(2 - \alpha C_f)$.

\end{lemma}
\emph{Proof:} See \cite{koshal11}, Lemma 4.4. \hfill $\blacksquare$

Next we establish an inequality that will be used to prove Theorem \ref{thm:muclose} below. 

\begin{lemma} \label{lem:coco}
Let Assumption \ref{as:functions} hold. Then the inequality 
\begin{equation*}
(\mu_2 - \mu_1)^T(-g(x_2) + g(x_1)) \geq \frac{\nu}{M_g^2}\|x_2 - x_1\|^2
\end{equation*}
is satisfied by all $x_1, x_2 \in X$ and $\mu_1, \mu_2 \in \R^m_{+}$,
where $M_g = \max_{x \in X} \|\nabla g(x)\|$. 
\end{lemma}
\emph{Proof:} See \cite{koshal11}, Lemma 4.1. \hfill $\blacksquare$

We now state one of the key theoretical results of the paper. 

\begin{theorem} \label{thm:muclose}
Let Assumptions \ref{as:functions}--\ref{as:gradglip} hold.
Suppose that each agent uses regularization parameter $\nu > 0$ and the cloud
uses regularization parameter $\epsilon > 0$. Let $\alpha$ be bounded as in Lemma \ref{lem:geoconv} and let the dual stepsize $\tau$ be bounded according to
\begin{equation*}
\tau < \min\left\{\frac{2\nu}{M_g^2 + 2\epsilon\nu}, \frac{2\epsilon}{1 + \epsilon^2}\right\}.
\end{equation*}
Define the constant
\begin{equation*}
q_d := (1 - \tau\epsilon)^2 + \tau^2,
\end{equation*}
which is in the set $(0,1)$ by the definition of $\tau$. 
Then for all $t \in \N$ we have
\begin{multline} \label{eq:thm1statement}
\|\mu(t+1) - \hmve\|^2 \leq q_d^{t+1}\|\mu(0) - \hmve\|^2 \\ + \sum_{i=1}^{t}
q_d^{i-1}\Big(q_dM_g^2M_x^2q_p^{d(t-i)} + 2\tau^2M_g^2M_x^2q_p^{d(t-i)/2}\Big),
\end{multline}
where $M_{x} = \max_{x, y \in X} \|x - y\|$ and $q_p$ is as defined in Lemma \ref{lem:geoconv}. 
\end{theorem}
\noindent \emph{Proof}: 
For economy of notation we define
\begin{equation*}
g^{t-1} = g(x^{t-1}(d(t-1))).
\end{equation*}
We see that 
\begin{equation} \label{eq:pdv1}
\mu(t+1) = \Pi_{\dv}\left[\mu(t) + \tau\big(g^{t-1} - \ep\mu(t)\big)\right]
\end{equation}
and
\begin{equation} \label{eq:pdv2}
\hmve = \Pi_{\dv}\left[\hmve + \tau\big(g(\hxve) - \ep\hmve\big)\right].
\end{equation}
The projection operator $\pmu[\cdot]$ is non-expansive and hence using Equations \eqref{eq:pdv1} and \eqref{eq:pdv2}
gives
\begin{multline} 
\|\mu(t+1) - \hmve\|^2 \leq \\ (1 - \tau\ep)^2\|\mu(t) - \hmve\|^2 + \tau^2\|g(\hxve) - g^{t-1}\|^2 \\
                       -2\tau(1 - \tau\ep)(\mu(t) - \hmve)^T(g(\hxve) - g^{t-1}) \label{eq:bigp3}.
\end{multline}

We now place a bound on the term $\tau^2\|g(\hxve) - g^{t-1}\|^2$. 
Defining $x^{t-1}_* = (x^{t-1,T}_{*,1}, \ldots, x^{t-1,T}_{*,N})^T$,
we see that
\begin{multline*}
\|g(\hxve) - g^{t-1}\|^2 = \|g(\hxve) - g(x^{t-1}_{*}) + g(x^{t-1}_{*}) - g^{t-1}\|^2 \\
          \leq \|g(\hxve) - g(x^{t-1}_*)\|^2 + \|g(x^{t-1}_*) - g^{t-1}\|^2 \\
          + 2\|g(\hxve) - g(x^{t-1}_*)\|\|g(x^{t-1}_*) - g^{t-1}\|.
\end{multline*}
From the Lipschitz property of $g$, we see that
\begin{equation*}
\|g(\hxve) - g(x^{t-1}_*)\| \leq M_g \|\hxve - x^{t-1}_*\|
\end{equation*}
and that
\begin{equation*}
\|g(x^{t-1}_*) - g^{t-1}\| \leq M_g \|x^{t-1}_* - x^{t-1}(d(t-1))\|,
\end{equation*}
which now together give
\begin{multline} \label{eq:term1bound}
\|g(\hxve) - g^{t-1}\|^2 \leq 
\|g(\hxve) - g(x^{t-1}_*)\|^2 + \|g(x^{t-1}_*) - g^{t-1}\|^2 \\
 + 2M_g^2\|x^{t-1}_* - x^{t-1}(d(t-1))\|\|\hxve - x^{t-1}_*\|.
\end{multline}

Examining the term $(\mu(t) - \hmve)^T(g(\hxve) - g^{t-1})$ in Equation \eqref{eq:bigp3}, we add and subtract
the same term inside the second set of parentheses, giving
\begin{multline}
(\mu(t) - \hmve)^T(g(\hxve) - g^{t-1}) \\
 = (\mu(t) - \hmve)^T\big(g(\hxve) - g(x^{t-1}_*) + g(x^{t-1}_*) - g^{t-1}\big) \\
 \geq \frac{\nu}{M_g^2}\|g(\hxve) - g(x^{t-1}_*)\|^2 + (\mu(t) - \hmve)^T(g(x^{t-1}_*) - g^{t-1}) \label{eq:uglycoco},
\end{multline}
where the last inequality follows from Lemma \ref{lem:coco}. 
Multiplying Equation \eqref{eq:uglycoco} by $-2\tau(1 - \tau\ep)$ then gives
\begin{multline}
-2\tau(1 - \tau\ep)(\mu(t) - \hmve)^T\big(g(\hxve) - g^{t-1}\big) \\
\leq  -2\tau(1 - \tau\ep)\frac{\nu}{M_g^2}\|g(\hxve) - g(x^{t-1}_*)\|^2 + \tau^2\|\mu(t) - \hmve\|^2  \\ + (1 - \tau\ep)^2\|g(x^{t-1}_*) - g^{t-1}\|^2 \label{eq:helper2}.
\end{multline}

Substituting Equations \eqref{eq:term1bound} and \eqref{eq:helper2} into Equation \eqref{eq:bigp3} and
combining like terms gives
\begin{multline} \label{eq:bigugly}
\|\mu(t+1) - \hmve\|^2 \leq \big((1 - \tau\ep)^2 + \tau^2\big)\|\mu(t) - \hmve\|^2 \\
+ \left(\tau^2 - 2\tau(1 - \tau\ep)\frac{\nu}{M_g^2}\right)\|g(\hxve) - g(x^{t-1}_*)\|^2
\\ + \big((1 - \tau\ep)^2 + \tau^2\big)\|g(x^{t-1}_*) - g(x^{t-1}(d(t-1))\|^2
\\ + 2\tau^2M_g^2\|\hxve - x^{t-1}_*\|\|x^{t-1}_* - x^{t-1}(d(t-1))\|. 
\end{multline}
By Lemma \ref{lem:geoconv}, we see that
\begin{equation} \label{eq:reftolem}
\|x^{t-1}(d(t-1)) - x^{t-1}_*\| \leq q_p^{d(t-1)/2}\|x^{t-1}(0) - x^{t-1}_*\|. 
\end{equation}
Also, by definition of $M_g$ and $M_x$ we have
\begin{equation} \label{eq:setup1}
\|g(\hxve) - g(x^{t-1}_*)\|^2 \leq M_g^2\|\hxve - x^{t-1}_*\|^2 \leq M_g^2M_x^2,
\end{equation}
while using $M_g$ and Lemma \ref{lem:geoconv} gives
\begin{multline} \label{eq:ugh}
\|g(x^{t-1}_*) - g^{t-1}\|^2 \leq \\ M_g^2q_p^{d(t-1)}\|x^{t-1}(0) - x^{t-1}_*\|^2 \leq M_g^2M_x^2q_p^{d(t-1)}.
\end{multline}
Along the lines of Equation \eqref{eq:reftolem}, we also see that 
\begin{equation} \label{eq:setup2}
\|\hxve - x^{t-1}_*\|\|x^{t-1}_* - x^{t-1}(d(t-1))\| \leq M_x^2q_p^{d(t-1)/2}.
\end{equation}

Using Equations \eqref{eq:reftolem}, \eqref{eq:setup1}, \eqref{eq:ugh}, and \eqref{eq:setup2} in
Equation \eqref{eq:bigugly} gives
\begin{multline} \label{eq:ugly2}
\|\mu(t+1) - \hmve\|^2 \leq \\ \big((1 - \tau\ep)^2 + \tau^2\big)\|\mu(t) - \hmve\|^2
 + \left(\tau^2 - 2\tau(1 - \tau\ep)\frac{\nu}{M_g^2}\right)M_g^2M_x^2 \\
+ \left((1 - \tau\ep)^2 + \tau^2\right)M_g^2M_x^2q_p^{d(t-1)} 
+ 2\tau^2M_g^2M_x^2q_p^{d(t-1)/2}.
\end{multline}
Using that
\begin{equation*}
\tau < \min\left\{\frac{2\nu}{M_g^2 + 2\ep\nu}, \frac{2\ep}{1 + \ep^2}\right\},
\end{equation*}
we have
\begin{equation*}
(1 - \tau\ep)^2 + \tau^2 < 1
\end{equation*}
and
\begin{equation*}
\tau^2 - 2\tau(1 - \tau\ep)\frac{\nu}{M_g^2} < 0.
\end{equation*}
Setting $q_d = (1 - \tau\ep)^2 + \tau^2$,  Equation \eqref{eq:ugly2} gives
\begin{multline} \label{eq:yay}
\|\mu(t+1) - \hmve\|^2 \leq \\ q_d\|\mu(t) - \hmve\|^2
+ q_dM_g^2M_x^2q_p^{d(t-1)}
+ 2\tau^2M_g^2M_x^2q_p^{d(t-1)/2}.
\end{multline}
Applying Equation \eqref{eq:yay} recursively then gives the desired result. \hfill $\blacksquare$

The sum in the statement of Theorem \ref{thm:muclose} contains terms that will eventually
become small by virtue of the leading term being $q_d^{i-1}$ and $i$ getting large. 
The interpretation of Theorem
\ref{thm:muclose} is that vectors $\mu(t)$ computed by the cloud will eventually
become close to $\hmve$, though the distance between them will never become zero. In fact, 
after a long time 
this distance
is dominated by the first few terms of the summation in Equation \eqref{eq:thm1statement} 
precisely because
the term containing $q_d^{t+1}$ goes to zero. We now bound the distance between primal
vectors and their optima.

\begin{theorem} \label{thm:xseq}
Let Assumptions 1-4 hold. 
Then for the sequence of primal vectors $\{x^t\big(d(t)\big)\}_{t \in \N}$ generated by
Algorithm \ref{alg:dist} we have
\begin{equation*}
\|x^t\big(d(t)\big) - \hxve\| \leq q_p^{d(t)/2}M_x + \frac{M_g}{\nu}\|\mu(t) - \hmve\|
\end{equation*}
along with
\begin{equation*} 
\max\{0, g_j(x^t\big(d(t)\big)\} \leq M_g \left(q_p^{d(t)/2}M_x + \frac{M_g}{\nu}\|\mu(t) - \hmve\|\right).
\end{equation*}
\end{theorem}

\noindent \emph{Proof}: Using Lemmas \ref{lem:geoconv} and \ref{lem:coco}, we see that
\begin{align}\label{eq:xbound}
\|x^{t}(d(t)) - \hxve\| &\leq \|x^{t}(d(t)) - x^t_*\| + \|x^t_* - \hxve\| \\
&\leq q_p^{d(t)/2}M_x + \frac{M_g}{\nu}\|\mu(t) - \hmve\|. 
\end{align}

Using the convexity of each constraint function $g_j$, we see that
\begin{align*}
g_j\big(x^t(d(t))\big) &\leq g_j(\hxve) + \grad g(x^t(d(t)))^T(x^t(d(t)) - \hxve) \\
        &\leq \|\grad g(\hxve)\|\|x^t(d(t)) - \hxve\| \\ 
        &\leq M_g\big(q_p^{d(t)/2}M_x + \frac{M_g}{\nu}\|\mu(t) - \hmve\|\big),
\end{align*}
for all $j \in J$, where we have used the bound on $\|x^t(d(t)) - \hxve\|$ from Equation~\eqref{eq:xbound}. 
\hfill $\blacksquare$

The first half of Theorem \ref{thm:xseq} says that $x_i^t$ will eventually become close to $\hxvei$, with the degree
of closeness determined in part by the distance between $\mu(t)$ and $\hmve$. The second half makes a similar statement
about the extent of any constraint violations, namely that the degree of any constraint violation depends upon
the distance from $\mu(t)$ to $\hmve$. Both statements in Theorem \ref{thm:xseq} say that
the length of delays between dual updates can be beneficial when it is long, 
 though naturally longer delays also
require more time for convergence. 

While the point $(\hxve, \hmve)$ is not necessarily 
a saddle point of the original (unregularized) Lagragian, 
 it is guaranteed to be sufficiently close
 when $\nu$ and $\epsilon$ are small enough; an extended discussion of this matter is in~\cite[Section 3.2]{koshal11}. 
In effect, Algorithm \ref{alg:dist} lets the agents approach $\hxve$ which itself it not far from an optimal solution to 
Problem~1 in terms of the optimal function value and small feasibility violation of the functional constraints.

%% file: 5-results.tex
\section{Experimental Results} \label{sec:numer}
\begin{figure}
\centering
\includegraphics[width=3.4in]{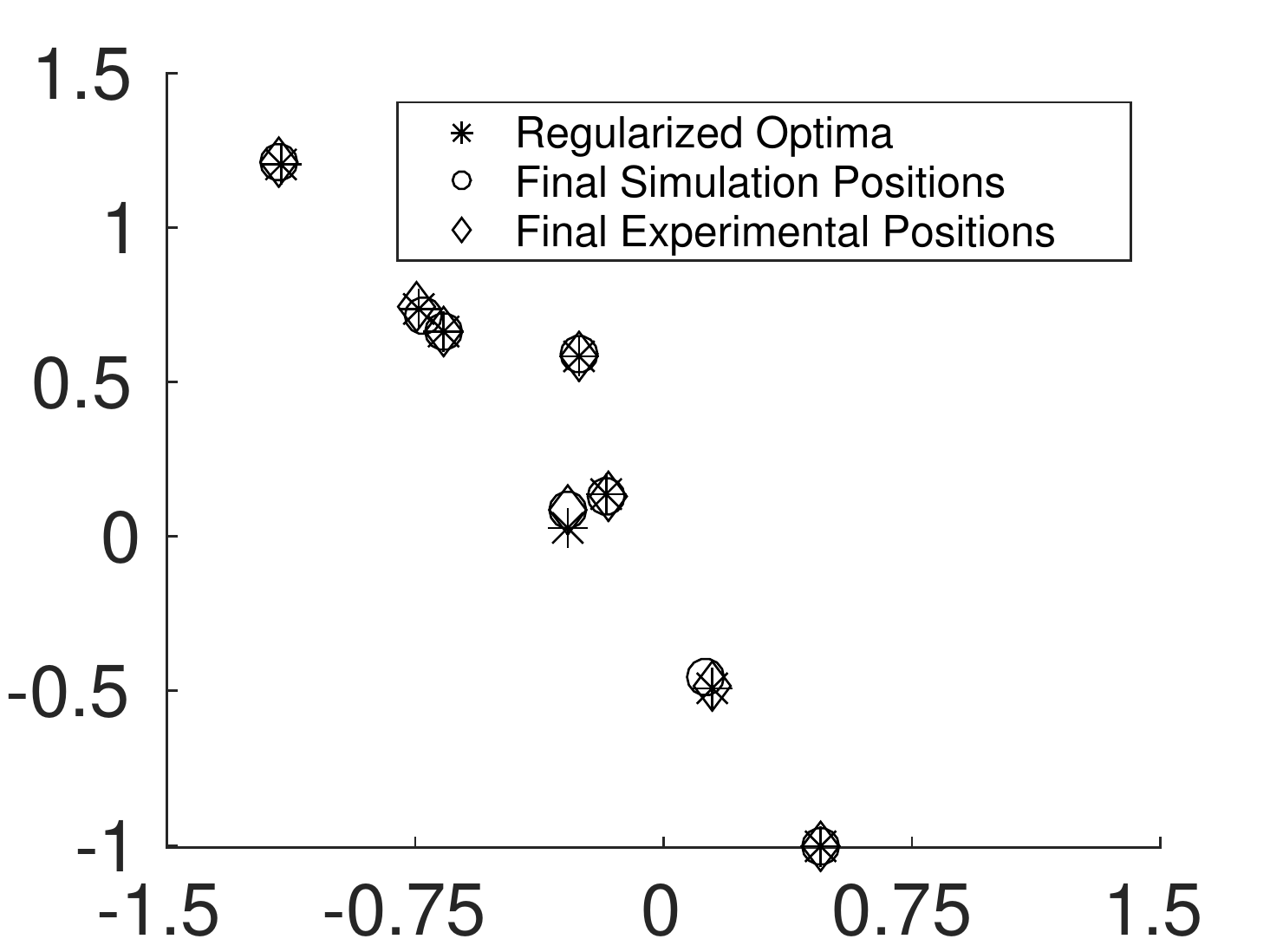}
\caption{A plot of $\hxve$, the final simulation positions, and the final robot positions, shown as asterisks, circles, and diamonds, respectively.}
\label{fig:finalpos}
\end{figure}
Algorithm \ref{alg:dist} was simulated for and then run with $8$ planar agents.
We outline the problem and cover the simulation results, and
then present the experimental results.

\begin{figure}
\centering
\includegraphics[width=3.4in,natwidth=1180,natheight=761]{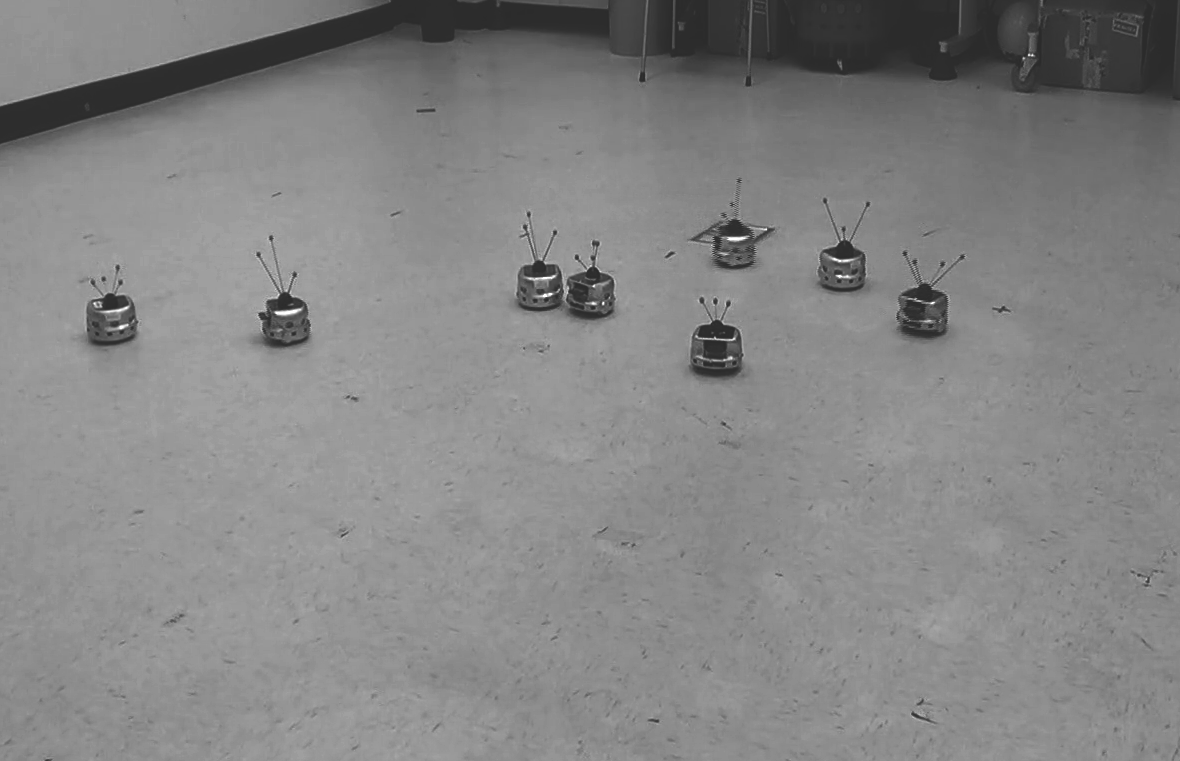}
\caption{Team of $8$ Khepera III robots executing the cloud-based algorithm.}
\label{fig:robots}
\end{figure}

All agents are planar so that $x_i \in \R^2$ for all $i$
and $x \in \R^{16}$. 
The sum of the per-agent objective functions is
\begin{multline*}
\sum_{i=1}^8 f_i(x_i) = \|x_1\|^2 + 
\left\|x_2 - \left(\ns\ns\begin{array}{r}-1 \\ 1\end{array}\ns\ns\right)\right\|^2 +
\left\|x_3 - \left(\ns\ns\begin{array}{r} 0.2  \\ -0.6\end{array}\ns\ns\right)\right\|^2 \\ +
\left\|x_4 - \left(\ns\ns\begin{array}{r}-1.4 \\ 1.4\end{array}\ns\ns\right)\right\|^2 +
\left\|x_5 - \left(\ns\ns\begin{array}{r}-0.1 \\ 0.5 \end{array}\ns\ns\right)\right\|^2 +
\left\|x_6 - \left(\ns\ns\begin{array}{r}-0.7  \\ 0.7 \end{array}\ns\ns\right)\right\|^2 \\ +
(x_{7,1} - 0.5)^2 + x_{7,2} - 1.1 +
(x_{8,1} + 0.3)^2 + x_{8,2}^4.
\end{multline*}
The non-separable term in the cost is 
\begin{equation*}
c(x) = \frac{1}{200}\big(\|x_1 - x_4\|^2 + \|x_1 - x_8\|^2 + \|x_4 - x_8\|^2\big).
\end{equation*}
As before the total cost used was $f(x) = \sum_{i=1}^8 f_i(x_i) + c(x)$. 
The constraints used were
\begin{equation*}
g(x) = \left(\begin{array}{l} \|x_1 - x_2\|^2 - 0.6 \\ \|x_1 - x_5\|^2 - 1.2 \\ \|x_7 - x_8\|^2 - 1.8 \\ \|x_1 - x_3\|^2 - 0.4 \\ \|x_4 - x_6\|^2 - 0.9 \end{array}\right) \leq 0,
\end{equation*}
and each agent was confined to the box $[-1.5,1.5]\times[-1.0,1.5]$, giving
\begin{equation*}
X = \prod_{i=1}^8 [-1.5,1.5]\times[-1.0,1.5].
\end{equation*}

The constants needed to solve this problem were computed (approximately) numerically and are shown in Table \ref{tab:constants}. 
\begin{table}[!htbp]
\centering
\begin{tabular}{|c|c|} \hline
Symbol & Value \\ \hline
$L_f$ & $26.9982$ \\ \hline
$L_g$ & $4.2426$ \\  \hline
$M_{\nu}$ & $8.7430$ \\ \hline
$M_g$ & $13.5277$ \\ \hline
$M_x$ & $\sqrt{122}$ \\ \hline
\end{tabular} 
\caption{Constants associated with Problem \ref{prob:vi}}
\label{tab:constants}
\end{table}
\vspace{-0.5cm}

The regularization parameters were chosen to be $\nu = \epsilon = 0.1$, giving $C_f = 64.191$. The primal and dual stepsizes
were chosen to be
\begin{equation*}
\alpha = 0.9\cdot\frac{2}{C_f} = 0.02804
\end{equation*}
and
\begin{equation*}
\tau = 0.9\cdot \min\left\{\frac{2\nu}{M_g^2 + 2\epsilon\nu}, \frac{2\epsilon}{1 + \epsilon^2}\right\} = 9.835\cdot 10^{-4}.
\end{equation*}

All delays had length determined by a random integer drawn from a uniform distribution on the integers betwen $10$ and $100$,
inclusive. Algorithm \ref{alg:dist} was run until the agents had computed $250,000$ state updates, during which time
the cloud compued $4,539$ dual updates. 

In simulation, the initial total distance between the agents' positions and their regularized optima, $\|x^0(0) - \hxve\|$, was
$2.6024$ and their final total distance, $\|x^{4539}(250000) - \hxve\|$, was $0.0843$.
In addition, after only $25,000$ iterations
the total distance of the agents to $\hxve$ was $0.1392$, indicating that fewer steps can be taken while still
achieving an acceptable ending state. 
In the dual space, the final distance to the regularized optimum was $\|\mu(4,539) - \hmve\| = 0.0393$,
indicating close convergence in the dual space as well. 

This problem was executed on a team of $8$ Khepera III 
robots; 
these robots are pictured in Figure \ref{fig:robots} where they are $14$ seconds into the experimental run. Position data was gathered using 
an OptiTrack 
motion capture system and the cloud-based algorithm was used to generate position waypoints for the agents.
The experiment was run
until the robots and cloud completed $75,000$ total updates; this point was reached in the middle of a communications cycle, and that cycle was
allowed to finish, giving $73,720$ total state updates by each agent and $1,340$ dual updates by the cloud. The final error in the primal
space was $\|x^{1340}(73720) - \hxve\| = 0.0689$ and the final error in the dual space was $\|\mu(1340) - \hmve\| = 0.0587$, indicating
close convergence of the robots to $(\hxve, \hmve)$. 
A plot of the regularized optima, simulation results, and experimental results is shown in Figure \ref{fig:finalpos},
indicating close agreement among the three sets of data plotted there.

%% file: 6-conclusion.tex
\section{Conclusion} \label{sec:concl}
We presented a hybrid centralized/decentralized algorithm
for solving multi-agent nonlinear programs with inequality constraints.
To do this, we used a Tikhonov regularization of the problem and a computing
regime that spread computations across the agents and a cloud computer. 
The architectural model incorporated communications delays in the system
and approximate convergence of the algorithm was proven. Experimental
results were provided to show the applicability of these results.

%% file: delayed-cloud.bbl
\begin{thebibliography}{10}

\bibitem{chiang07}
Mung Chiang, S.H. Low, A.R. Calderbank, and J.C. Doyle.
\newblock Layering as optimization decomposition: A mathematical theory of
  network architectures.
\newblock {\em Proceedings of the IEEE}, 95(1):255--312, Jan 2007.

\bibitem{cortes05}
Jorge Cort\'{e}s, Sonia Martínez, and Francesco Bullo.
\newblock Spatially-distributed coverage optimization and control with
  limited-range interactions.
\newblock {\em ESAIM: Control, Optimisation and Calculus of Variations},
  11:691--719, 10 2005.

\bibitem{gharesifard14}
B.~Gharesifard and J.~Cortes.
\newblock Distributed continuous-time convex optimization on weight-balanced
  digraphs.
\newblock {\em Automatic Control, IEEE Transactions on}, 59(3):781--786, March
  2014.

\bibitem{guo02}
Y.~Guo and L.E. Parker.
\newblock A distributed and optimal motion planning approach for multiple
  mobile robots.
\newblock In {\em Robotics and Automation, 2002. Proceedings. ICRA '02. IEEE
  International Conference on}, volume~3, pages 2612--2619, 2002.

\bibitem{hale14}
M.T. Hale and M.~Egerstedt.
\newblock Cloud-based optimization: A quasi-decentralized approach to
  multi-agent coordination.
\newblock In {\em Decision and Control, 2014. Proceedings of the 53rd IEEE
  Conference on}. IEEE, 2014.
\newblock To appear.

\bibitem{nazari14}
M~Honarvar~Nazari, Zak Costello, Mohammad~Javad Feizollahi, Santiago Grijalva,
  and Magnus Egerstedt.
\newblock Distributed frequency control of prosumer-based electric energy
  systems.
\newblock {\em Power Systems, IEEE Transactions on}, 29, November 2014.

\bibitem{kelly98}
F.~Kelly, A.~Maulloo, and D.~Tan.
\newblock Rate control in communication networks: shadow prices, proportional
  fairness and stability.
\newblock In {\em Journal of the Operational Research Society}, volume~49,
  1998.

\bibitem{khan09}
M.~Khan, G.~Pandurangan, and V.S.A. Kumar.
\newblock Distributed algorithms for constructing approximate minimum spanning
  trees in wireless sensor networks.
\newblock {\em Parallel and Distributed Systems, IEEE Transactions on},
  20(1):124--139, Jan 2009.

\bibitem{konnov07}
Igor Konnov.
\newblock {\em Equilibrium models and variational inequalities}, volume 210.
\newblock Elsevier, 2007.

\bibitem{koshal11}
Jayash Koshal, Angelia Nedi\'{c}, and Uday~V. Shanbhag.
\newblock Multiuser optimization: Distributed algorithms and error analysis.
\newblock {\em SIAM Journal on Optimization}, 21(3):1046--1081, 2011.

\bibitem{kuhn51}
H.~W. Kuhn and A.~W. Tucker.
\newblock Nonlinear programming.
\newblock In {\em Proceedings of the Second Berkeley Symposium on Mathematical
  Statistics and Probability}, pages 481--492, Berkeley, Calif., 1951.
  University of California Press.

\bibitem{lobel11}
I.~Lobel and A.~Ozdaglar.
\newblock Distributed subgradient methods for convex optimization over random
  networks.
\newblock {\em Automatic Control, IEEE Transactions on}, 56(6):1291--1306, June
  2011.

\bibitem{mitra94}
Debasis Mitra.
\newblock An asynchronous distributed algorithm for power control in cellular
  radio systems.
\newblock In {\em Wireless and Mobile Communications}, pages 177--186.
  Springer, 1994.

\bibitem{nedic09}
A.~Nedi\'{c} and A.~Ozdaglar.
\newblock Distributed subgradient methods for multi-agent optimization.
\newblock {\em Automatic Control, IEEE Transactions on}, 54(1):48--61, Jan
  2009.

\bibitem{notarstefano09}
G.~Notarstefano and F.~Bullo.
\newblock Network abstract linear programming with application to cooperative
  target localization.
\newblock In {\em Modelling, Estimation and Control of Networked Complex
  Systems}, Understanding Complex Systems, pages 177--190. 2009.

\bibitem{panait05}
Liviu Panait and Sean Luke.
\newblock Cooperative multi-agent learning: The state of the art.
\newblock {\em Autonomous Agents and Multi-Agent Systems}, 11(3):387--434,
  2005.

\bibitem{rabbat04}
M.~Rabbat and R.~Nowak.
\newblock Distributed optimization in sensor networks.
\newblock In {\em Information Processing in Sensor Networks, 2004. IPSN 2004.
  Third International Symposium on}, pages 20--27, April 2004.

\bibitem{soltero13}
Daniel~E Soltero, Mac Schwager, and Daniela Rus.
\newblock Decentralized path planning for coverage tasks using gradient descent
  adaptive control.
\newblock {\em The International Journal of Robotics Research}, 2013.

\bibitem{srivastava11}
K.~Srivastava and A.~Nedic.
\newblock Distributed asynchronous constrained stochastic optimization.
\newblock {\em Selected Topics in Signal Processing, IEEE Journal of},
  5(4):772--790, Aug 2011.

\bibitem{trigoni12}
Niki Trigoni and Bhaskar Krishnamachari.
\newblock {Sensor network algorithms and applications: Introduction}.
\newblock {\em {Philosophical Transactions of the Royal Scoeity A -
  Mathematical, Physical, and Engineering Sciences}}, {370}({1958,
  SI}):{5--10}, {JAN 13} {2012}.

\bibitem{tsitsiklis84}
John Tsitsiklis.
\newblock {\em {Problems in Decentralized Decision making and Computation}}.
\newblock PhD thesis, Massachusetts Institute of Technology, 1984.

\bibitem{uzawa58}
H.~Uzawa.
\newblock Iterative methods in concave programming.
\newblock {\em {Studies in Linear and Non-Linear Programming}}, 1958.

\bibitem{wei13}
Ermin Wei, Asuman Ozdaglar, and Ali Jadbabaie.
\newblock A distributed newton method for network utility maximization i:
  Algorithm.
\newblock {\em Automatic Control, IEEE Transactions on}, 58(9):2162--2175,
  2013.

\end{thebibliography}
